# Hankel determinants of q-exponential polynomials

*Johann Cigler*

Fakultät für Mathematik
Universität Wien
A-1090 Wien, Nordbergstraße 15

johann.cigler@univie.ac.at

**Abstract**
We give simple proofs for the Hankel determinants of $q$ – exponential polynomials.

Let $S(n,k)$ be the Stirling numbers of the second kind. Christian Radoux ([6]) has shown that the Hankel determinants of the exponential polynomials $B_n(x) = \sum_{k=0}^{n} S(n,k) x^k$ are given by

$$\det\left(B_{i+j}(x)\right)_{i,j=0}^{n-1} = x^{\binom{n}{2}} \prod_{j=0}^{n-1} j!. \tag{1}$$

In [2] I have proved some $q$ – analogues of this result. Then Richard Ehrenborg [4] has given a combinatorial proof of one of these $q$ – analogues.
In this paper I want to show that these $q$ – analogues in some sense have simpler proofs than the original case.

We use the usual notations: For $n \in \mathbb{N}$ let $[n] = \dfrac{1-q^n}{1-q}$. The $q$ – factorial is the product $[1] \cdot [2] \cdots [n]$ and the $q$ – binomial coefficient $\begin{bmatrix} n \\ k \end{bmatrix}$ is defined by $\begin{bmatrix} n \\ k \end{bmatrix} = \dfrac{[n]!}{[k]![n-k]!}$ for $0 \le k \le n$ and $\begin{bmatrix} n \\ k \end{bmatrix} = 0$ else.

The $q$ – Stirling numbers $S[n,k]$ of the second kind are defined by

$$S[n,k] = S[n-1,k-1] + [k]S[n-1,k] \tag{2}$$

with $S[0,k] = [k=0]$ and $S[n,0] = [n=0]$.

There are two natural $q$-analogues of the exponential polynomials:

$$\varphi_n(x) = \sum_{k=0}^{n} S[n,k] x^k \tag{3}$$

and

$$\Phi_n(x) = \sum_{k=0}^{n} q^{\binom{k}{2}} S[n,k] x^k. \tag{4}$$

Our aim is a simple proof of the following theorems.

**Theorem 1**
The Hankel determinants of the $q$-exponential polynomials $\varphi_n(x)$ are given by

$$\det\left(\varphi_{i+j}(x)\right)_{i,j=0}^{n-1} = x^{\binom{n}{2}} q^{\binom{n}{3}} \prod_{j=0}^{n-1} [j]! \tag{5}$$

and

$$\det\left(\varphi_{i+j+1}(x)\right)_{i,j=0}^{n-1} = x^{\binom{n+1}{2}} q^{\binom{n+1}{3}} \prod_{j=0}^{n-1} [j]!. \tag{6}$$

**Theorem 2**
The Hankel determinants of the $q$-exponential polynomials $\Phi_n(x)$ are

$$\det\left(\Phi_{i+j}(x)\right)_{i,j=0}^{n-1} = q^{2\binom{n}{3}} x^{\binom{n}{2}} \prod_{j=0}^{n-1} \left([j]!((1-q)x;q)_j\right) \tag{7}$$

and

$$\det\left(\Phi_{i+j+1}(x)\right)_{i,j=0}^{n-1} = q^{2\binom{n+1}{3}} x^{\binom{n+1}{2}} \prod_{j=0}^{n-1} \left([j]!((1-q)x;q)_j\right). \tag{8}$$

Here $(x;q)_n$ is defined by $(x;q)_n = \prod_{j=0}^{n-1}\left(1 - q^j x\right)$.

The key for the simpler proofs is the well-known identity

$$(q-1)^{n-k} S[n,k] = \sum_i (-1)^{n-i} \binom{n}{i} \begin{bmatrix} i \\ k \end{bmatrix}. \tag{9}$$

The equivalence of (2) and (9) may be seen from the following computation:

$$\sum_i (-1)^{n-i}\binom{n}{i}\begin{bmatrix}i\\k\end{bmatrix} = \sum_i (-1)^{n-i}\left(\binom{n-1}{i-1}+\binom{n-1}{i}\right)\begin{bmatrix}i\\k\end{bmatrix}$$

$$= \sum_i (-1)^{n-i}\binom{n-1}{i-1}\begin{bmatrix}i\\k\end{bmatrix} + \sum_i (-1)^{n-i}\binom{n-1}{i}\begin{bmatrix}i\\k\end{bmatrix}$$

$$= \sum_i (-1)^{n-1-i}\binom{n-1}{i}\left(\begin{bmatrix}i+1\\k\end{bmatrix}-\begin{bmatrix}i\\k\end{bmatrix}\right) = \sum_i (-1)^{n-1-i}\binom{n-1}{i}q^{i+1-k}\begin{bmatrix}i\\k-1\end{bmatrix}$$

$$= \sum_i (-1)^{n-1-i}\binom{n-1}{i}\begin{bmatrix}i\\k-1\end{bmatrix} + \sum_i (-1)^{n-1-i}\binom{n-1}{i}\begin{bmatrix}i\\k-1\end{bmatrix}(q^{i+1-k}-1)$$

$$= \sum_i (-1)^{n-1-i}\binom{n-1}{i}\begin{bmatrix}i\\k-1\end{bmatrix} + \sum_i (-1)^{n-1-i}\binom{n-1}{i}\begin{bmatrix}i\\k\end{bmatrix}(q^k-1)$$

We need some other well-known results:

**Lemma 1**
*For given sequences $s(n)$ and $t(n)$ define $a(n,k)$ by*

$$\begin{aligned}a(0,k) &= [k=0] \\ a(n,0) &= s(0)a(n-1,0)+t(0)a(n-1,1) \\ a(n,k) &= a(n-1,k-1)+s(k)a(n-1,k)+t(k)a(n-1,k+1).\end{aligned} \quad (10)$$

*Then the Hankel determinant $\det(a(i+j,0))_{i,j=0}^{n-1}$ is given by*

$$\det(a(i+j,0))_{i,j=0}^{n-1} = \prod_{i=1}^{n-1}\prod_{k=0}^{i-1} t(k). \quad (11)$$

For a proof see e.g. [3].

**Remark**
In most cases we start with $a(n)=a(n,0)$ and want to find the corresponding $s(n)$ and $t(n)$. It is then convenient to compute the first values of the orthogonal polynomials $p(n,x)$ (cf. [3] (1.10)) and their Favard resolution [3] (1.11) and try to guess $s(n)$ and $t(n)$. Then guess the explicit form of $a(n,k)$. Afterwards it remains to verify (10) in order to obtain a rigorous proof.

**Lemma 2**
*Define the binomial transform of a sequence $(a_n)$ by $BIN(a_n)=(b_n)$ with $b_n = \sum_{k=0}^{n}\binom{n}{k}a_k$.*

*Then*

$$\det(a_{i+j}) = \det(b_{i+j}). \quad (12)$$

A simple proof can be found in [7].

Further observe that
$$\det\left(x^{i+j}a_{i+j}\right)_{i,j=0}^{n-1} = x^{n(n-1)}\det\left(a_{i+j}\right)_{i,j=0}^{n-1}. \tag{13}$$

From (9) we get
$$\varphi_n(x) = \sum_{k=0}^{n} S[n,k]x^k = \sum_k x^k(q-1)^{k-n}\sum_i(-1)^{n-i}\binom{n}{i}\begin{bmatrix}i\\k\end{bmatrix}$$
$$= \sum_i(-1)^{n-i}\binom{n}{i}(q-1)^{-n}\sum_k(q-1)^k x^k\begin{bmatrix}i\\k\end{bmatrix}.$$

In terms of the Rogers-Szegö polynomials $r_n(x)$ defined by
$$r_n(x) = \sum_{k=0}^{n}\begin{bmatrix}n\\k\end{bmatrix}x^k \tag{14}$$

this means
$$\varphi_n(x) = \frac{1}{(1-q)^n}\sum_i\binom{n}{i}(-1)^i r_i((q-1)x). \tag{15}$$

By (13) and (12) this implies that
$$\det\left(\varphi_{i+j}(x)\right)_{i,j=0}^{n-1} = \frac{1}{(q-1)^{\binom{n}{2}}}\det\left(r_{i+j}(x)\right)_{i,j=0}^{n-1}. \tag{16}$$

Therefore we have only to determine the Hankel determinants of the Rogers-Szegö polynomials. These are also well-known (cf. [5]), but can also be obtained in a trivial way from (11).
The Rogers-Szegö polynomials satisfy the recurrence (cf. e.g. [1] )
$$r_n(x) = (x+1)r_{n-1}(x) + (q^{n-1}-1)xr_{n-2}(x). \tag{17}$$

Let now $s(k) = q^k(x+1)$ und $t(k) = q^k x(q^{k+1}-1)$.
Then it is easily verified that the corresponding $a(n,k)$ are given by
$$a(n,k) = \begin{bmatrix}n\\k\end{bmatrix}r_{n-k}(x). \tag{18}$$

We have only to check that (10) holds:
$$\begin{bmatrix}n\\k\end{bmatrix}r_{n-k}(x) = \begin{bmatrix}n-1\\k-1\end{bmatrix}r_{n-k}(x) + q^k(x+1)\begin{bmatrix}n-1\\k\end{bmatrix}r_{n-k-1}(x) + q^k x(q^{k+1}-1)\begin{bmatrix}n-1\\k+1\end{bmatrix}r_{n-k-2}(x)$$
or
$$\begin{bmatrix}n-1\\k\end{bmatrix}r_{n-k}(x) = (x+1)\begin{bmatrix}n-1\\k\end{bmatrix}r_{n-k-1}(x) + x(q^{n-k-1}-1)\begin{bmatrix}n-1\\k\end{bmatrix}r_{n-k-2}(x)$$

which is immediate from (17).

Therefore by (11) we get

**Lemma 3**

$$\det\left(r_{i+j}(x)\right)_{i,j=0}^{n-1} = x^{\binom{n}{2}} q^{\binom{n}{3}} (q-1)^{\binom{n}{2}} \prod_{j=0}^{n-1} [j]!. \tag{19}$$

This immediately implies (5).

Let now $D$ denote the $q$-differentiation operator, defined by $Df(x) = \dfrac{f(x) - f(qx)}{(1-q)x}$.

Then

$$\varphi_n(x) = x(1+D)\varphi_{n-1}(x). \tag{20}$$

For
$$\varphi_n(x) = \sum_k S[n,k]x^k = \sum_k S[n-1,k-1]x^k + \sum_k S[n-1,k][k]x^k$$
$$= x\varphi_{n-1}(x) + xD\varphi_{n-1}(x).$$

This implies

$$\varphi_n(x) = \left(x(1+D)\right)^n 1. \tag{21}$$

Let now $\varepsilon$ be the linear operator defined by $\varepsilon f(x) = f(qx)$. Then

$$\varepsilon x D \varepsilon^{-1} x^n = \varepsilon x D \frac{x^n}{q^n} = [n]x^n = xDx^n, \text{ d.h. } \varepsilon x D \varepsilon^{-1} = xD.$$

Therefore we get
$$\varphi_{n+1}(x) = \left(x(1+D)\right)^{n+1} 1 = x\left((1+D)x\right)^n (1+D)1 = x\left((1+D)x\right)^n 1.$$

Now $Dxx^n = [n+1]x^n = (q[n]+1)x^n = (qxD+1)x^n$ for all $n$. Thus $Dx = 1 + qxD$.

This implies
$$\varphi_{n+1}(x) = x\left((1+D)x\right)^n 1 = x\left(1 + x + qxD\right)^n 1 = x\varepsilon^{-1}\left(1 + qx(1+D)\right)^n 1$$

and therefore the recurrence relation

$$\varphi_{n+1}(x) = x \sum_{k=0}^{n} \binom{n}{k} q^k \varphi_k\left(\frac{x}{q}\right). \tag{22}$$

By (12) this implies (6). Therefore Theorem 1 is proved.

Formula (5) can be slightly generalized in the following form:
Let $c \geq 1$ be an integer and

$$\psi_n(x,c) = \left(x^c + xD\right)^n 1. \tag{23}$$

Then we get

$$\psi_n(x,c) = \sum_{k=0}^{n} S[n,k,q^c] x^{kc} \left([c]\right)^{n-k}, \tag{24}$$

where $S[n,k,q^c] = S[n,k]_{q \to q^c}$.

For
$$\left(x^c + xD\right)\sum_{k=0}^{n-1} S[n-1,k,q^c]x^{kc}[c]^{n-k-1} =$$

$$\sum_k S[n-1,k,q^c]x^{(k+1)c}\left([c]\right)^{n-k-1} + \sum_k S[n-1,k,q^c]x^{kc}[kc]\left([c]\right)^{n-k-1}$$

$$= \sum_k x^{kc}\left([c]\right)^{n-k}\left(S[n-1,k-1,q^c] + \frac{[kc]}{[c]}S[n-1,k,q^c]\right) = \sum_{k=0}^n S[n,k,q^c]x^{kc}\left([c]\right)^{n-k}.$$

From

$$\psi_n(x,c) = \sum_{k=0}^n S[n,k,q^c]x^{kc}\left([c]\right)^{n-k} = \sum_k x^{ck}(q-1)^{k-n}\sum_i (-1)^{n-i}\binom{n}{i}\begin{bmatrix}i\\k\end{bmatrix}_c$$

$$= \sum_i (-1)^{n-i}\binom{n}{i}(q-1)^{-n}\sum_k (q-1)^k x^{ck}\begin{bmatrix}i\\k\end{bmatrix}_c$$

we get in the same way as above

$$\det\left(\psi_{i+j}(x,c)\right)_{i,j=0}^{n-1} = \frac{1}{(q-1)^{2\binom{n}{2}}}\det\left(r_{i+j}((q-1)x^c)\big|_{q^c}\right)_{i,j=0}^{n-1} = \frac{(q^c-1)^{\binom{n}{2}}}{(q-1)^{\binom{n}{2}}}x^{c\binom{n}{2}}q^{c\binom{n}{3}}\prod_{j=0}^{n-1}[j]_c!$$

This gives

$$\det\left(\psi_{i+j}(x,c)\right)_{i,j=0}^{n-1} = [c]^{\binom{n}{2}}x^{c\binom{n}{2}}q^{c\binom{n}{3}}\prod_{j=0}^{n-1}[j]_c! \tag{25}$$

In the same way as above we get

$$\psi_{n+1}(x,c) = x^c\sum_{k=0}^n \binom{n}{k}q^{kc}\psi_k\left(\frac{x}{q},c\right)[c]^{n-k}. \tag{26}$$

This implies

$$\det\left(\psi_{i+j+1}(x,c)\right) = x^{nc}q^{c\binom{n}{2}}\det\left(\psi_{i+j}(x,c)\right). \tag{27}$$

**Proof of Theorem 2**

Consider now the $q$ – exponential polynomials $\Phi_n(x)$.

From (9) we get
$$\Phi_n(x) = \sum_{k=0}^{n} q^{\binom{k}{2}} S[n,k] x^k = \sum_k q^{\binom{k}{2}} x^k (q-1)^{k-n} \sum_i (-1)^{n-i} \binom{n}{i} \begin{bmatrix} i \\ k \end{bmatrix}$$
$$= \sum_i (-1)^{n-i} \binom{n}{i} (q-1)^{-n} \sum_k q^{\binom{k}{2}} (q-1)^k x^k \begin{bmatrix} i \\ k \end{bmatrix}.$$

Observing that $(x;q)_n = \prod_{j=0}^{n-1}(1-q^j x) = \sum_{k=0}^{n}(-1)^k \begin{bmatrix} n \\ k \end{bmatrix} q^{\binom{k}{2}} x^k$ we see that

$$\Phi_n(x) = \frac{1}{(1-q)^n} \sum_i \binom{n}{i} (-1)^i ((1-q)x;q)_i. \tag{28}$$

This implies that
$$\det \left( \Phi_{i+j}(x) \right)_{i,j=0}^{n-1} = \frac{1}{(1-q)^{2\binom{n}{2}}} \det \left( ((1-q)x;q)_{i+j} \right)_{i,j=0}^{n-1}. \tag{29}$$

Therefore we have only to determine the Hankel determinants of the polynomials $(x;q)_n$.

Let $s(n) = q^n + q^{n-1} x \left(1 - q^n (1+q)\right)$
and $t(n) = q^{2n}(1-q^{n+1})x(1-q^n x)$.
If we define $a(n,k)$ by (10) then we get
$$a(n,k) = \begin{bmatrix} n \\ k \end{bmatrix} (q^k x; q)_{n-k}.$$

To prove this assertion we have to verify that
$$\begin{bmatrix} n \\ k \end{bmatrix} (q^k x; q)_{n-k} = \begin{bmatrix} n-1 \\ k-1 \end{bmatrix} (q^{k-1} x; q)_{n-k} + s(k) \begin{bmatrix} n-1 \\ k \end{bmatrix} (q^k x; q)_{n-k-1} + t(k) \begin{bmatrix} n-1 \\ k+1 \end{bmatrix} (q^{k+1} x; q)_{n-k-2}$$

holds. This is equivalent with
$$\begin{bmatrix} n \\ k \end{bmatrix} (q^k x; q)_{n-k} = \begin{bmatrix} n-1 \\ k-1 \end{bmatrix} (q^{k-1} x; q)_{n-k} + \left( q^k + q^{k-1} x - q^{2k-1} x - q^{2k} x \right) \begin{bmatrix} n-1 \\ k \end{bmatrix} (q^k x; q)_{n-k-1}$$
$$+ (q^{2k} - q^{n+k-1}) x \begin{bmatrix} n-1 \\ k \end{bmatrix} (q^k x; q)_{n-k-1}$$

or
$$\begin{bmatrix} n \\ k \end{bmatrix} (q^k x; q)_{n-k} = \begin{bmatrix} n-1 \\ k-1 \end{bmatrix} (q^{k-1} x; q)_{n-k} + \left( q^k + q^{k-1} x - q^{2k-1} x - q^{n+k-1} x \right) \begin{bmatrix} n-1 \\ k \end{bmatrix} (q^k x; q)_{n-k-1}$$
$$= \begin{bmatrix} n-1 \\ k-1 \end{bmatrix} (q^{k-1} x; q)_{n-k} + q^k \begin{bmatrix} n-1 \\ k \end{bmatrix} (q^{k-1} x; q)_{n-k} + q^{k-1} x (1-q^n) \begin{bmatrix} n-1 \\ k \end{bmatrix} (q^k x; q)_{n-k-1}$$
$$= \begin{bmatrix} n \\ k \end{bmatrix} (q^{k-1} x; q)_{n-k} + q^{k-1} x (1-q^{n-k}) \begin{bmatrix} n \\ k \end{bmatrix} (q^k x; q)_{n-k-1}$$

which is obviously true.

Thus we get

**Lemma 4**

$$\det\left((x;q)_{i+j}\right)_{i,j=0}^{n-1} = q^{2\binom{n}{3}}(1-q)^{\binom{n}{2}} x^{\binom{n}{2}} \prod_{j=0}^{n-1} [j]!\, (x;q)_j. \tag{30}$$

*and*

$$\det\left((x;q)_{i+j+1}\right)_{i,j=0}^{n-1} = q^{2\binom{n}{3}+\binom{n}{2}}(1-q)^{\binom{n}{2}} x^{\binom{n}{2}} \prod_{j=0}^{n-1} [j]!\, (x;q)_{j+1}. \tag{31}$$

The second equation follows from $(x;q)_{n+1} = (1-x)(qx;q)_n$.

(30) implies immediately (7).

**Remark**

The special case $x = q$ gives the well-known result

$$\det\left([i+j]!\right)_{i,j=0}^{n-1} = q^{\frac{n(n-1)(2n-1)}{6}} \prod_{j=0}^{n-1} [j]!^2. \tag{32}$$

As another special case we consider the $q$-Stirling numbers $s[n,k]$ of the first kind. They satisfy $s[n+1,k] = s[n,k-1] - [n]s[n,k]$ with initial values $s[0,k] = [k=0]$ and $s[n,0] = [n=0]$.
Then it is easily verified that

$$\sum_{k=0}^{n} s[n,k] x^k = \langle x \rangle_n := \prod_{j=0}^{n-1} (x - [j]). \tag{33}$$

From $\langle x \rangle_n = \dfrac{(1+(q-1)x)^n}{(q-1)^n} \left(\dfrac{1}{1+(q-1)x}; q\right)_n$

we get

$$\det\left(\langle x \rangle_{i+j}\right)_{i,j=0}^{n-1} = (-1)^{\binom{n}{2}} q^{2\binom{n}{3}} \prod_{j=0}^{n-1} [j]!\, \langle x \rangle_j. \tag{34}$$

In order to obtain the second Hankel determinant we observe that

$$\Phi_n(x) = \sum_k q^{\binom{k}{2}} S[n,k] x^k = \sum_k q^{\binom{k}{2}} S[n-1,k-1] x^k + \sum_k q^{\binom{k}{2}} S[n-1,k][k] x^k$$
$$= x\Phi_{n-1}(qx) + xD\Phi_{n-1}(x)$$

or

$$\Phi_n(x) = (x\varepsilon + xD)\Phi_{n-1}(x) = (x\varepsilon + xD)^n 1. \tag{35}$$

This gives

$$\Phi_{n+1}(x) = (x\varepsilon + xD)^{n+1} 1 = x((\varepsilon + D)x)^n (\varepsilon + D) 1 = x((\varepsilon + D)x)^n 1$$
$$= x(\varepsilon x + qxD + 1)^n 1 = x(q(x\varepsilon + xD) + 1)^n 1.$$

Therefore

$$\Phi_{n+1}(x) = x \sum_{k=0}^n \binom{n}{k} q^k \Phi_k(x). \tag{36}$$

By (12) this implies (8). Therefore Theorem 2 is proved.

This theorem can be generalized in the same way as above.
Let $c \geq 1$ be an integer and

$$\Psi_n(x,c) = (x^c \varepsilon + xD)^n 1. \tag{37}$$

Then we get

$$\Psi_n(x,c) = \sum_{k=0}^n S[n,k,q^c] x^{kc} q^{c\binom{k}{2}} ([c])^{n-k}, \tag{38}$$

where $S[n,k,q^c] = S[n,k]_{q \to q^c}$.

From

$$\Psi_n(x,c) = \sum_{k=0}^n S[n,k,q^c] x^{kc} q^{c\binom{k}{2}} ([c])^{n-k} = \sum_k x^{ck} q^{c\binom{k}{2}} (q-1)^{k-n} \sum_i (-1)^{n-i} \binom{n}{i} \begin{bmatrix} i \\ k \end{bmatrix}_c$$
$$= \sum_i (-1)^{n-i} \binom{n}{i} (q-1)^{-n} \sum_k (q-1)^k x^{ck} q^{c\binom{k}{2}} \begin{bmatrix} i \\ k \end{bmatrix}_c$$

we get in the same way as above

$$\det \left( \Psi_{i+j}(x,c) \right)_{i,j=0}^{n-1} = \frac{1}{(q-1)^{2\binom{n}{2}}} \det \left( ((1-q)x^c; q^c)_{i+j} \right)_{i,j=0}^{n-1} = [c]^{\binom{n}{2}} x^{c\binom{n}{2}} q^{2c\binom{n}{3}} \prod_{j=0}^{n-1} [j]_c! ((1-q)x^c; q^c)_j$$

This gives

$$\det\left(\Psi_{i+j}(x,c)\right)_{i,j=0}^{n-1} = [c]^{\binom{n}{2}} x^{c\binom{n}{2}} q^{2c\binom{n}{3}} \prod_{j=0}^{n-1} [j]_c!\left((1-q)x^c; q^c\right)_j. \tag{39}$$

Furthermore we get

$$\Psi_{n+1}(x,c) = x^c \sum_{k=0}^{n} \binom{n}{k} [c]^{n-k} q^{kc} \Psi_k(x,c). \tag{40}$$

This implies finally

$$\det\left(\Psi_{i+j+1}(x,c)\right)_{i,j=0}^{n-1} = x^{nc} q^{2c\binom{n}{2}} \det\left(\Psi_{i+j}(x,c)\right)_{i,j=0}^{n-1}. \tag{41}$$

**Remark**

The method used for the Rogers-Szegö polynomials can also be applied to the $q$-analogue of the Hermite polynomials defined by

$$H_n(x) = xH_{n-1}(x) - [n-1]H_{n-2}(x) \tag{42}$$

with initial values
$H_0(x) = 1$ and $H_1(x) = x.$

If we choose $s(n) = q^n x$ and $t(n) = -q^n[n+1]$ we get

$$a(n,k) = \begin{bmatrix} n \\ k \end{bmatrix} H_{n-k}(x).$$

We have only to check that

$$\begin{bmatrix} n \\ k \end{bmatrix} H_{n-k}(x) = \begin{bmatrix} n-1 \\ k-1 \end{bmatrix} H_{n-k}(x) + s(k) \begin{bmatrix} n-1 \\ k \end{bmatrix} H_{n-k-1}(x) + t(k) \begin{bmatrix} n-1 \\ k+1 \end{bmatrix} H_{n-k-2}(x)$$

or

$$q^k \begin{bmatrix} n-1 \\ k \end{bmatrix} H_{n-k}(x) = q^k x \begin{bmatrix} n-1 \\ k \end{bmatrix} H_{n-k-1}(x) - q^k[k+1] \begin{bmatrix} n-1 \\ k+1 \end{bmatrix} H_{n-k-2}(x)$$

$$= q^k x \begin{bmatrix} n-1 \\ k \end{bmatrix} H_{n-k-1}(x) - q^k[n-k-1] \begin{bmatrix} n-1 \\ k \end{bmatrix} H_{n-k-2}(x).$$

But this is true because of (42).

Therefore we get

$$\det\left(H_{i+j}(x)\right)_{i,j=0}^{n-1} = (-1)^{\binom{n}{2}} q^{\binom{n}{3}} \prod_{j=0}^{n-1} [j]!. \tag{43}$$

# References

[1]  J. Cigler, Elementare q-Identitäten, Sém. Lotharingien Comb. B05 A 1982

[2]  J. Cigler, Eine Charakterisierung der q- Exponentialpolynome, Sitzber. OeAW , 208 (1999), 143-157,  http://hw.oeaw.ac.at/sitzungsberichte_und_anzeiger_collection

[3]  J. Cigler, Hankel determinants of Schröder-like numbers, arXiv:0901.4680

[4]  R. Ehrenborg, Determinants involving q-Stirling numbers, Adv. in Appl. Math. 31 (2003), 630 - 642

[5]  Q.-H. Hou , A. Lascoux and Y.-P. Mu, Continued fractions for Rogers-Szegö polynomials, Numerical Algorithms 35 (2004), 81-90

[6]  C. Radoux,  Calcul effectif de certains déterminants de Hankel, Bull. Soc. Math. Belg., Sér. B, 31 (1979), 49-55

[7]  D. Zeilberger, An umbral approach to the Hankel transform of sequences, 2005 http://www.math.rutgers.edu/~zeilberg/mamarim/mamarimPDF/hankel.pdf